\documentclass{amsart}[12]

\usepackage[margin=1.0in]{geometry}
\ignorespacesafterend
\usepackage{amsmath,amsthm,amssymb,array,enumerate,parskip,graphicx,etoolbox,tabularx}
\usepackage[all]{xy}

\usepackage{xcolor} 
\usepackage{tikz-cd}

\usepackage[pagebackref]{hyperref}
\definecolor{dark-red}{rgb}{0.5,0.15,0.15}
\definecolor{dark-blue}{rgb}{0.15,0.15,0.6}
\definecolor{dark-green}{rgb}{0.15,0.6,0.15}
\hypersetup{
    colorlinks, linkcolor=dark-red,
    citecolor=dark-blue, urlcolor=dark-green
}

\usepackage[capitalise]{cleveref}
\usepackage{comment}

\newcommand{\Z}{\mathbb{Z}} 
\newcommand{\Q}{\mathbb{Q}} 

\newcommand \ZZ {{\mathbb Z}}

\usepackage{breqn}

 \theoremstyle{definition}

\title[]{A note on the procedure to find the generic polynomial of a quotient (closely following Adelmann)}
\author{Rachel Davis}
\address{University of Wisconsin-Madison}
\email{rachel.davis@wisc.edu}

\thanks{ }

\begin{document}

\maketitle

\section*{Introduction} 
There are 3 examples in these notes. The first one is the standard example of the cubic resolvent of a quartic. The second example is exactly from Adelmann \cite{Adelmann} and gives a defining polynomial corresponding to the unique $S_4$-quotient of $\mathrm{GL}_2(\Z/4\Z)$. The splitting field of the Adelmann polynomial over $\Q$ is a subfield of the 4-division field of an elliptic curve, that contains the 2-division field of the elliptic curve. The third example is new and needed in the study of the field theory of quaternion origami. Associated to an elliptic curve defined over $\Q$, with a rational point, is a degree 8 polynomial whose Galois group is a subgroup of $\mathrm{Hol}(Q_8)$. Three defining polynomials corresponding to the three $S_4$-quotients of  $\mathrm{Hol}(Q_8)$ are given.

\section{Warm-up example: Resolvent cubic}

\subsection{Roots}

Suppose we are given a quartic polynomial $f(x)=x^4+a_3x^3+a_2x^2+a_1x+a_0$ with roots $\alpha_1, \alpha_2, \alpha_3, \alpha_4$. Suppose also, that the Galois group of the splitting field of $f(x)$ is $S_4$. Group-theoretically, there is exactly 1 $S_3$-quotient of $S_4$. One can associate to $f(x)$, a resolvent cubic, i.e. a degree 3 polynomial $g(x)$ such that the splitting field of $g(x)$ is contained in the splitting field of $f(x)$ and the Galois group of $g(x)$ is $S_3$. There are formulas readily available for this polynomial, but to help systematize this process, we will follow the full procedure used in Adelmann \cite{Adelmann}. Let $G=S_4$, $H=V_4 \leftrightarrow \left\{ (1,4)(2,3), (1,3)(2,4)\right\}$, a normal subgroup, and $N=G/H$ as an abstract group isomorphic to $S_3$. We are looking for a subgroup $F$ of $G$ such that its projection $\overline{F}$ into $N$ has exactly 3 conjugates in $N$. We conclude that the normalizer of $\overline{F}$ in $N$ needs to have index 3. Since $N$ is isomorphic to $S_3$, the normalizer described above is isomorphic to $\Z/2\Z$. There are suitable $F$ having sizes 2, 4, and 8 respectively. In order for the resolvent polynomial to have degree 3, we require $|F|=8$. Since all of the 3 possible subgroups $F$ lead to the same resolvent polynomials, we may select any of these groups as our $F$, for example, take $F=\langle (3,4), (1,4)(2,3), (1,3)(2,4) \rangle $. A system of representatives for $G/F$ is given by $\left\{ (), (2,3), (2,4) \right\}$.   

Finally, we need a polynomial $P \in k[x_1,x_2, x_3,x_4]$ having stabilizer subgroup $\mathrm{Stab}_G(P)=F$ which should be homogenous and of lowest possible degree. There is no suitable polynomial in degree 1, but in degree 2, we find $P=\overline{R}_F(x_1x_2)=4(x_1x_2+x_3x_4)$. Then, we can take $\beta_1=\alpha_1\alpha_2+\alpha_3\alpha_4$, $\beta_2=\alpha_1\alpha_3+\alpha_2\alpha_4$, and $\beta_3=\alpha_1\alpha_4+\alpha_2\alpha_3$.

\subsection{Vieta's formulas}

By using Vieta's formulas for relating the roots and the coefficients of a polynomial, we can find the coefficients of the resolvent cubic $g(x)=( x-\beta_1 )( x-\beta_2 )( x-\beta_3 )=x^3 + (-\alpha_1\alpha_2 - \alpha_1\alpha_3 - \alpha_1\alpha_4 - \alpha_2\alpha_3 - \alpha_2\alpha_4 - \alpha_3\alpha_4)x^2 + (\alpha_1^2\alpha_2\alpha_3 +
    \alpha_1^2\alpha_2\alpha_4 + \alpha_1^2\alpha_3\alpha_4 + \alpha_1\alpha_2^2\alpha_3 + \alpha_1\alpha_2^2\alpha_4 + \alpha_1\alpha_2\alpha_3^2 + \alpha_1\alpha_2\alpha_4^2
    + \alpha_1\alpha_3^2\alpha_4 + \alpha_1\alpha_3\alpha_4^2 + \alpha_2^2\alpha_3\alpha_4 + \alpha_2\alpha_3^2\alpha_4 + \alpha_2\alpha_3\alpha_4^2)x -
    \alpha_1^3\alpha_2\alpha_3\alpha_4 - \alpha_1^2\alpha_2^2\alpha_3^2 - \alpha_1^2\alpha_2^2\alpha_4^2 - \alpha_1^2\alpha_3^2\alpha_4^2 -
    \alpha_1\alpha_2^3\alpha_3\alpha_4 - \alpha_1\alpha_2\alpha_3^3\alpha_4 - \alpha_1\alpha_2\alpha_3\alpha_4^3 - \alpha_2^2\alpha_3^2\alpha_4^2$.
    
For example, the $x^2$-coefficient of $fgx)$ is $-a_2$, the negative of the $x^2$-coefficient of $f$. The $x$-coefficient of $g$ is $a_1a_3 - 4a_0$. The constant coefficient of $g$ is $-a_1^2 + 4a_0a_2 - a_0a_3^2$. Therefore, the resolvent cubic is $g(x)=x^3-a_2x^2+(a_1a_3 - 4a_0)x+(-a_1^2 + 4a_0a_2 - a_0a_3^2)$. This agrees with a standard definition of a resolvent cubic polynomial of $f(x)$.   

\section{Adelmann example} 
\subsection{Roots}
This is verbatim from p. 104 of \cite{Adelmann}. Let $E: y^2=x^3+ax+b$ be an elliptic curve defined over $\Q$, the field of rational numbers and consider the degree 6 polynomial $A(Y)=\frac{1}{2}\Lambda_4(Y)=Y^6+5aY^4+20bY^3-5a^2Y^2-4abY-a^3-8b^2$. The Galois group of the splitting field of this polynomial over $\Q$ is a subgroup of $\mathrm{PGL}_2(\Z/4\Z)=\mathrm{GL}_2(\Z/4\Z)/\langle \pm I \rangle$. There is a unique $S_4$-quotient of this group. We would like to give the generic polynomial cutting out the splitting field of this quotient. Let $G=\mathrm{PGL}_2(\Z/4\Z)=\mathrm{GL}_2(\Z/4\Z)/\langle \pm I \rangle$. Let $H = \left\{  \left(\begin{array}{cc}
 1&  0\\
  0&  1

\end{array} \right),   \left(\begin{array}{cc}
 1&  2 \\
  2 & -1

\end{array} \right) \right\} \longleftrightarrow \left\{ (), (1,6)(2,4)(3,5) \right\}$. The subgroup $H$ of $G$ is a normal subgroup, the factor group $N=G/H$, as an abstract group, is isomorphic to $S_4$. Hence, the polynomial to be determined should have degree 4. Consequently, we are looking for a subgroup $F$ of $G$, so that its projection $\overline{F}$ into $N$ has exactly 4 conjugates within $N$. We conclude that the normalizer of $\overline{F}$ in $N$ needs to have index 4 in $N$. Since $N$ is isomorphic to $S_4$, the normalizer described is isomorphic to $S_3$, and we discover suitable $F$ having sizes 6 and 12, respectively. In order for the resolvent polynomial to have degree 4 we require $|F|=12$. Since all the 4 possible subgroups $F$ lead to the same resolvent polynomials, we may select any of these groups as our $F$, for example

$F= \left\{ \right.$
 \begin{equation*}
\begin{split}
& \left(\begin{array}{cc}
 1&  0\\
  0&  1 \end{array} \right),  
  \left(\begin{array}{cc}
 0&  -1\\
  1&  0 \end{array} \right),  \left(\begin{array}{cc}
 -1&  0\\
  1&  1 \end{array} \right),   \\
  \noindent & \left(\begin{array}{cc}
  -1 & 1 \\
  0 & 1 \end{array} \right) \left(\begin{array}{cc}
 1&  1\\
  1&  0
\end{array} \right),   \left(\begin{array}{cc}
 0&  1\\
  1&  -1 \end{array} \right), \\
  &  \left(\begin{array}{cc}
 1&  2\\
  2&  -1 \end{array} \right),  \left(\begin{array}{cc}
 2&  1\\
  1&  2 \end{array} \right), 
 \left(\begin{array}{cc}
 1&  2\\
  1&  -1 \end{array} \right),  \\
  & \left(\begin{array}{cc}
 1&  1\\
  2&  -1 \end{array} \right),  
\left(\begin{array}{cc}
 -1&  1\\
  1&  2 \end{array} \right),  \left(\begin{array}{cc}
 2&  1\\
  1&  1 \end{array} \right)  
\end{split}
\end{equation*}
$\left. \right\}$

 $F$ is generated by
 $$ \left\{  \left(\begin{array}{cc}
 1&  1\\
  1&  0 \end{array} \right),  \left(\begin{array}{cc}
 0&  -1\\
  1&  0 \end{array} \right) \right\} \longleftrightarrow \left\{ (1,2,5,6,4,3), (1,4)(2,6) \right\}. $$
  
  A system of representatives of $G/F$ is given by

  $$ \left\{  \left(\begin{array}{cc}
 1&  0\\
  0&  1 \end{array} \right),  \left(\begin{array}{cc}
 1&  0\\
  2&  -1 \end{array} \right),  \left(\begin{array}{cc}
 1&  2\\
  0&  -1 \end{array} \right),  \left(\begin{array}{cc}
 1&  0\\
  0&  -1 \end{array} \right) \right\} \longleftrightarrow \left\{ (), (1,6), (2,4), (3,5) \right\}. $$

Finally, we need a polynomial $P \in k[x_1, \ldots, x_6]$ having stabilizer subgroup $\mathrm{Stab}_G(P)=F$ which should be homogeneous and of lowest possible degree. There is no suitable polynomial in degree 1, but in degree 2 we find $P=\overline{R}_F(x_1x_2)$.

$P=\overline{R}_F(x_1x_2)=x_1x_2+x_2x_5+x_5x_6+x_4x_6+x_3x_4+x_1x_3+x_4x_6+x_5x_6+x_2x_5+x_1x_2+x_1x_3+x_3x_4=2(x_1x_2+x_2x_5+x_5x_6+x_4x_6+x_3x_4+x_1x_3)$.

Take $\beta_1=x_1x_2+x_2x_5+x_5x_6+x_4x_6+x_3x_4+x_1x_3$, $\beta_2=x_2x_6+x_2x_5+x_1x_5+x_1x_4+x_3x_4+x_3x_6$,
$\beta_3=x_1x_4+x_4x_5+x_5x_6+x_2x_6+x_2x_3+x_1x_3$,
$\beta_4=x_1x_2+x_2x_3+x_3x_6+x_4x_6+x_4x_5+x_1x_5$.

\subsection{Vieta's formulas}
By using Vieta's formulas relating the roots and the coefficients of a polynomial, we can find the coefficients of the resolvent quartic $g(x)=( x-\beta_1 )( x-\beta_2 )( x-\beta_3 )(x-\beta_4)$. This time, the formulas involve a little more than the symmetric polynomials, so there is a bit of additional work to do. Let $\omega_f$ be the evaluation homomorphism from $k[x_1, \ldots, x_6] \rightarrow k_f$ given by $x_i \mapsto \alpha_i$ where the $\alpha_i$s are the roots of $f$ (In this case. $f=A(Y)=\frac{1}{2}\Lambda_4(Y)$. 

We know that $\omega_f(\overline{R}_G(x_2x_4))$ and $\omega_f(\overline{R}_G(x_1x_2))$ sum to $5a$ since $\overline{R}_G(x_2x_4)$ and $\overline{R}_G(x_1x_2)$ sum to an elementary symmetric polynomial. 

Each of these elements is a fundamental invariant and so by Proposition 6.3.3 of Adelmann, \cite{Adelmann}, the coefficients of the resolvent polynomials $\mathcal{RP}_{S_6}(u)(Y)$ are symmetric polynomials in the $x_i$s and can therefore be expressed by elementary symmetric polynomials. By proposition 6.3.8 of Adelmann, the polynomial $\mathcal{R}_{S_6,u}(Y)$ is a polynomial whose coefficients are polynomials in the coefficients of the defining polynomial of $K$ and $\omega_f(u)$ appears as a simple linear factor. Adelmann gives $\mathcal{R}_{u_2}(Y)$, $\mathcal{R}_{u_3}(Y)$, $\mathcal{R}_{u_4}(Y)$ and each has a simple linear factor with root $\omega_f(u)$. The image of $\omega_f$ on all of the fundamental invariants can be determined by combining this information along with the images under $\omega_f$ of the elementary symmetric polynomials. As on p. 102 of Adelmann, \cite{Adelmann}, the result is the following:

\renewcommand\arraystretch{2.0}%
\begin{center}
  \begin{tabular}{| c | c| }
    \hline
    Fundamental Invariant & Image under $\omega_f$ \\ \hline 
   $u_1=\overline{R}_G(x_1)=6R_G(x_1)$ & $0$ \\ 
   $u_2=\overline{R}_G(x_2x_4)=3R_G(x_2x_4)$ & $a$ \\    
   $v_2=\overline{R}_G(x_1x_2)=12R_G(x_1x_2)$ & $4a$ \\  
    $u_3=\overline{R}_G(x_1x_2x_3)=8R_G(x_1x_2x_3)$ & $-8b$ \\  
    $u_4=\overline{R}_G(x_2x_2x_3x_4x_5)=3R_G(x_2x_3x_4x_5)$ & $-a^2$ \\
     $u_6=\overline{R}_G(x_1x_2x_2x_3x_4x_5x_6)=R_G(x_1x_2x_3x_4x_5x_6)$ & $-(a^3+8b^2)$ \\
        $w_3=\overline{R}_G(x_1x_2x_4)=12R_G(x_1x_2x_4)$ & $-12b$ \\ 
            $w_4=\overline{R}_G(x_1x_2x_3x_4)=12R_G(x_1x_2x_3x_4)$ & $-4a^2$ \\ 
                $w_5=\overline{R}_G(x_1x_2x_3x_4x_5)=6R_G(x_1x_2x_3x_4x_5)$ & $4ab$ \\ \hline
  \end{tabular}
\end{center}
\vspace{1cm}

\noindent Using these values gives that $\mathcal{R}_{G,P}(Y)=Y^4-8aY^3+24a^2Y^2+(224a^3+1728b^2)Y+272a^4+1728ab^2$ (See p. 105 of \cite{Adelmann}.). After performing the linear substitution $Y \mapsto Y+2a$ and introducing the quantity $\Delta=-16(4a^3+27b^2)$ (the discriminant of the underlying elliptic curve) we obtain the defining polynomial of the extension:

$$B(Y)=Y^4-4\Delta Y-12a\Delta.$$

\noindent Adelmann also gives the defining polynomials for the division field $\Q(E[n])/\Q$ using resultants. Adelmann \cite{Adelmann} provides these polynomials as follows:

$$T_n(X)=\mathrm{Res}_Y(\Gamma_n(Y), X^2-(Y^3-aY-b))$$

\noindent where $$\Gamma_n(X)=\displaystyle \prod_{d|n} A_d(X,Y)^{\mu(n/d)}$$
{\setlength{\abovedisplayskip}{0pt}%
 \begin{flalign*}
A_1&=1\\
 A_2&=2y\\
 A_3&=3x^4+6ax^2+12bx-a^2\\
  A_4&=4y(x^6+5ax^4+20bx^3-5a^2x^2-4abx-8b^2-a^3)\\
   2yA_{2m}&=A_m(A_{m+2}A_{m-1}^2-A_{m-2}A_{m+1}^2)\\ 
   A_{2m+1}&=A_{m+2}A_m^3-A_{m-1}A_{m+1}^3
   \end{flalign*}
  
\noindent Example: $n=4$
$$T_4(X)=\mathrm{Res}_Y(\Gamma_4(Y),X^2-(Y^3-aY-b))$$
\noindent where $\Gamma_4(X) = A_2(X,Y)^{-1}A_4(X,Y)=2(X^6+5aX^4+20bX^3-5a^2X^2-4abX-8b^2-a^3)$.

\noindent Then 
\begin{gather*}
\begin{split}
T_4(X)=X^{12} + 54bX^{10} + (132a^3 + 891b^2)X^8 + (432a^3b + 2916b^3)X^6 +\\
        (-528a^6 - 7128a^3b^2 - 24057b^4)X^4 + (864a^6b + 11664a^3b^3 +
        39366b^5)X^2\\ 
        - 64a^9 - 1296a^6b^2 - 8748a^3b^4 - 19683b^6
        \end{split}
        \end{gather*}

\section{$\mathrm{Hol}(Q_8)$ has 3 $S_4$-quotients }

\subsection{Roots}

Let $E: y^2=x^3+ax+b$ be an elliptic curve defined over $\Q$ with $P=(z,w) \in E(\Q)$ and consider the degree 8 polynomial $f(x)=f_{E,P, Q_8}(x)=x^8-8wx^6+6(2az+3b)x^4-(4a^3+27b^2)$. The Galois group of the splitting field of this polynomial over $\Q$ is a subgroup of $\mathrm{Hol}(Q_8)$ \cite{DavisGoins}. There are 3 distinct $S_4$-quotients of this group. We would like to give the generic polynomials cutting out the splitting fields corresponding to each of these quotients. Let $G=\mathrm{Hol}(Q_8)$. 

Then, $G$ is $\mathrm{TransitiveGroup}(8,40)$, as a permutation group. Write $G$ as \\
$\mathrm{sub}\langle S_8 | (1,3,4,8,7,5)(2,6), (2,6)(3,7), (2,4,7)(3,6,5), (3,7)(4,5), (1,7,4,2,8,3,5,6) \rangle$.

Then, consider the following 3 subgroups of $G$:

 $H_1=\langle (1,2,8,6)(3,5,7,4), (1,7,8,3)(2,5,6,4), (1,8)(2,6)(3,7)(4,5) \rangle$, and \\
 $H_2 = \langle  (1,4,8,5)(2,3,6,7), (1,2,8,6)(3,4,7,5), (1,8)(2,6)(3,7)(4,5)\rangle$, and \\
   $H_3= \langle (2,6)(4,5), (2,6)(3,7), (1,8)(2,6)(3,7)(4,5) \rangle$. \\

 Then, for $1 \leq i \leq 3$, $G/H_i \simeq S_4$. Here,  and $H_1 \simeq H_2 \simeq Q_8$ and $H_3 \simeq (\ZZ/2\ZZ)^3$.
 
 Let $N_i=G/H_i$. Then, $N_i=S_4$ for $i=1,2,3$. We are looking for degree 48 subgroups $F_i$ such that the normalizers of $\overline{F}_i$ in $N_i$ have exactly 4 conjugates. We will take

     $F_1=\langle (2,4)(3,7)(5,6), (2,3,5)(4,6,7), (1,2,8,6)(3,5,7,4),(1,7,8,3)(2,5,6,4), (1,8)(2,6)(3,7)(4,5)
 \rangle$
 
 $F_2 = \langle (2,5)(3,7)(4,6), (2,7,4)(3,5,6), (1,4,8,5)(2,3,6,7), (1,2,8,6)(3,4,7,5), (1,8)(2,6)(3,7)(4,5) \rangle$,
    
     $F_3=\langle(2,4,6,5), (2,7,4)(3,5,6), (2,6)(4,5), (2,6)(3,7), (1,8)(2,6)(3,7)(4,5) \rangle$
    
    One can take representatives for the cosets $G/F_i$ as follows (note that there are many different ways to choose a set of transversals for the $F_i$s):

        $$ \langle (), (3, 7)(4, 5), (1, 2, 5, 3, 8, 6, 4, 7),
    (1, 3, 5, 8, 7, 4)(2, 6) \rangle, $$
            $$ \langle (), (1,3,4,8,7,5)(2,6), (2,6)(3,7), (1,5,7,8,4,3)(2,6) \rangle, $$ and
           $$ \langle (), (1, 3, 4, 8, 7, 5)(2, 6), (1, 5, 7, 8, 4, 3)(2, 6),(1, 6, 5, 3, 8, 2, 4, 7) \rangle $$
    
    respectively. 
    
Finally, we need polynomials $P_i \in k[x_1, \ldots, x_8]$ having stabilizer subgroup $\mathrm{Stab}_G(P_i)=F_i$ which should be homogeneous and of lowest possible degree. We find $P_3=x_1x_8$ by studying the primary invariants of $F_3$.  Then, we can take $\gamma_1=x_1x_8$, $\gamma_2=x_2x_6$, $\gamma_3=x_3x_7$, $\gamma_4=x_4x_5$.

For $P_1$ and $P_2$, the primary invariants of $F_1$ and $F_2$ suggest using

$P_1=R_{F_1}(x_1x_2x_3)=x_1x_2x_3 + x_1x_2x_4 + x_1x_2x_5 + x_1x_3x_5 + x_1x_3x_6 + x_1x_4x_6 + x_1x_4x_7 + x_1x_5x_7 + x_1x_6x_7 + x_2x_3x_4 + x_2x_3x_8 + 
    x_2x_4x_7 + x_2x_5x_7 + x_2x_5x_8 + x_2x_7x_8 + x_3x_4x_6 + x_3x_4x_8 + x_3x_5x_6 + x_3x_5x_8 + x_4x_6x_8 + x_4x_7x_8 + x_5x_6x_7 + x_5x_6x_8 + x_6x_7x_8$ and

$P_2=R_{F_2}(x_1x_2x_3)=x_1x_2x_3 + x_1x_2x_4 + x_1x_2x_7 + x_1x_3x_5 + x_1x_3x_6 + x_1x_4x_6 + x_1x_4x_7 + x_1x_5x_6 + x_1x_5x_7 + x_2x_3x_4 + x_2x_3x_5 
        + x_2x_4x_8 + x_2x_5x_7 + x_2x_5x_8 + x_2x_7x_8 + x_3x_4x_6 + x_3x_4x_8 + x_3x_5x_8 + x_3x_6x_8 + x_4x_6x_7 + x_4x_7x_8 + 
        x_5x_6x_7 + x_5x_6x_8 + x_6x_7x_8$.
        
  Let $\alpha_i$ be the $G$-conjugates of $P_1$ with each $x_i$s evaluated with the $r_i$s (the roots of $f(x)$. \\
  
  $\alpha_1=r_1r_2r_3 + r_1r_2r_4 + r_1r_2r_5 + r_1r_3r_5 + r_1r_3r_6 + r_1r_4r_6 + r_1r_4r_7 + r_1r_5r_7 + r_1r_6r_7 + r_2r_3r_4 + r_2r_3r_8  + r_2r_4r_7 + r_2r_5r_7 + r_2r_5r_8 + r_2r_7r_8 + r_3r_4r_6 + r_3r_4r_8 + r_3r_5r_6 + r_3r_5r_8 + r_4r_6r_8 + r_4r_7r_8+ r_5r_6r_7 + r_5r_6r_8 + r_6r_7r_8$,\\
   $\alpha_2= r_1r_2r_4 + r_1r_2r_5 + r_1r_2r_7 + r_1r_3r_4 + r_1r_3r_5 + r_1r_3r_6 + r_1r_4r_7 + r_1r_5r_6 + r_1r_6r_7 + r_2r_3r_4 + r_2r_3r_5 + r_2r_3r_8 + r_2r_4r_8 + r_2r_5r_7 + r_2r_7r_8 + r_3r_4r_6 + r_3r_5r_8 + r_3r_6r_8 + r_4r_6r_7 + r_4r_6r_8 + r_4r_7r_8 + r_5r_6r_7 + r_5r_6r_8 + r_6r_7r_8$, \\ $
   \alpha_3= r_1r_2r_3 + r_1r_2r_4 + r_1r_2r_7 + r_1r_3r_4 + r_1r_3r_5 + r_1r_4r_6 + r_1r_5r_6 + r_1r_5r_7 + r_1r_6r_7 + r_2r_3r_5 + r_2r_3r_8 + r_2r_4r_7 + r_2r_4r_8 + r_2r_5r_7 + r_2r_5r_8 + r_3r_4r_6 + r_3r_4r_8 + r_3r_5r_6 + r_3r_6r_8 + r_4r_6r_7 + r_4r_6r_8 + r_5r_6r_7 + r_5r_7r_8 + r_6r_7r_8$,\\
$\alpha_4= r_1r_2r_3 + r_1r_2r_5 + r_1r_2r_7 + r_1r_3r_4 + r_1r_3r_6 + r_1r_4r_6 + r_1r_4r_7 + r_1r_5r_6 + r_1r_5r_7 + r_2r_3r_4 + r_2r_3r_5 + r_2r_4r_7 + r_2r_4r_8 + r_2r_5r_8 + r_2r_7r_8 +r_3r_4r_8 + r_3r_5r_6 + r_3r_5r_8 + r_3r_6r_8 + r_4r_6r_7 + r_4r_6r_8 +r_5r_6r_7 + r_5r_7r_8 + r_6r_7r_8$.\\
         
  Let $\beta_i$ be the $G$-conjugates of $P_2$ with each $x_i$s evaluated with the $r_i$s (the roots of $f(x)$).  \\
  $\beta_1=r_1r_2r_3 + r_1r_2r_4 + r_1r_2r_7 + r_1r_3r_5 + r_1r_3r_6 + r_1r_4r_6 + r_1r_4r_7 + r_1r_5r_6 + r_1r_5r_7 + r_2r_3r_4 + r_2r_3r_5 + r_2r_4r_8 + r_2r_5r_7 + r_2r_5r_8 + r_2r_7r_8 + r_3r_4r_6 + r_3r_4r_8 + r_3r_5r_8 + r_3r_6r_8 + r_4r_6r_7 + r_4r_7r_8 + r_5r_6r_7 + r_5r_6r_8 + r_6r_7r_8$,\\
 $\beta_2= r_1r_2r_3 + r_1r_2r_5 + r_1r_2r_7 + r_1r_3r_4 + r_1r_3r_5 + r_1r_4r_6 + r_1r_4r_7 + r_1r_5r_6 + r_1r_6r_7 + r_2r_3r_4 + r_2r_3r_8 + r_2r_4r_7 + r_2r_4r_8 + r_2r_5r_7 + r_2r_5r_8 + r_3r_4r_6 + r_3r_5r_6 + r_3r_5r_8 + r_3r_6r_8 + r_4r_6r_8 + r_4r_7r_8 + 
        r_5r_6r_7 + r_5r_7r_8 + r_6r_7r_8$, \\
         $\beta_3= r_1r_2r_4 + r_1r_2r_5 + r_1r_2r_7 + r_1r_3r_4 + r_1r_3r_5 + r_1r_3r_6 + r_1r_4r_6 + r_1r_5r_7 + r_1r_6r_7 + r_2r_3r_4 + r_2r_3r_5  + r_2r_3r_8 + r_2r_4r_7 + r_2r_5r_8 + r_2r_7r_8 + r_3r_4r_8 + r_3r_5r_6 + r_3r_6r_8 + r_4r_6r_7 + r_4r_6r_8 + r_4r_7r_8 +  r_5r_6r_7 + r_5r_6r_8 + r_5r_7r_8$, \\
         $\beta_4= r_1r_2r_3 + r_1r_2r_4 + r_1r_2r_5 + r_1r_3r_4 + r_1r_3r_6 + r_1r_4r_7 + r_1r_5r_6 + r_1r_5r_7 + r_1r_6r_7 + r_2r_3r_5 + r_2r_3r_8  + r_2r_4r_7 + r_2r_4r_8 + r_2r_5r_7 + r_2r_7r_8 + r_3r_4r_6 + r_3r_4r_8 + r_3r_5r_6 + r_3r_5r_8 + r_4r_6r_7 + r_4r_6r_8 + p1
        r_5r_6r_8 + r_5r_7r_8 + r_6r_7r_8$.
        
 The problem with these roots is that they are zero. For example, $\alpha_1=(x_1+x_8)(x_2x_3+x_2x_5+x_4x_7)+(x_2+x_6)(x_1x_4+x_5+x_7+x_7x_8)+(x_3+x_7)(x_1x_5+x_1x_6+x_2x_4)+(x_4+x_5)(x_3x_6+x_3x_8+x_6x_8)$ and since $x_1+x_8=0$, $x_2+x_6=0$, $x_3+x_7=0$, and $x_4+x_5=0$, $\alpha_1=0$ and similarly, $\alpha_2=0$, $\alpha_3=0$, $\alpha_4$. Similarly, the $\beta_i$ corresponding to the conjugates of $P_2$ evaluated at the roots of $f(x)$ are also zero.

Instead, we will use the following secondary invariants of $F_2$ and $F_3$, respectively (\cite{magma}) to construct roots:
 
     $P_1=x_1^2x_2x_7 + x_1^2x_3x_4 + x_1^2x_5x_6 + x_1x_2^2x_7 + x_1x_2x_7^2 + x_1x_3^2x_4 + x_1x_3x_4^2 + x_1x_5^2x_6 + 
    x_1x_5x_6^2 + x_2^2x_3x_5 + x_2^2x_4x_8 + x_2x_3^2x_5 + x_2x_3x_5^2 + x_2x_4^2x_8 + x_2x_4x_8^2 + x_3^2x_6x_8 + 
    x_3x_6^2x_8 + x_3x_6x_8^2 + x_4^2x_6x_7 + x_4x_6^2x_7 + x_4x_6x_7^2 + x_5^2x_7x_8 + x_5x_7^2x_8 + x_5x_7x_8^2$
    and    
   $P_2=x_1^2x_2x_3 + x_1^2x_4x_6 + x_1^2x_5x_7 + x_1x_2^2x_7 + x_1x_2x_4^2 + x_1x_3^2x_6 + x_1x_3x_5^2 + x_1x_4x_7^2 + 
    x_1x_5x_6^2 + x_2^2x_3x_5 + x_2^2x_4x_8 + x_2x_3^2x_4 + x_2x_5^2x_7 + x_2x_5x_8^2 + x_2x_7^2x_8 + x_3^2x_5x_8 + 
    x_3x_4^2x_6 + x_3x_4x_8^2 + x_3x_6^2x_8 + x_4^2x_7x_8 + x_4x_6^2x_7 + x_5^2x_6x_8 + x_5x_6x_7^2 + x_6x_7x_8^2$

 Then, the conjugates of these elements evaluated at the roots of $f(x)$ are the following:
 
  $\alpha_1=r_1^2r_2r_7 + r_1^2r_3r_4 + r_1^2r_5r_6 + r_1r_2^2r_7 + r_1r_2r_7^2 + r_1r_3^2r_4 
        + r_1r_3r_4^2 + r_1r_5^2r_6 + r_1r_5r_6^2 + r_2^2r_3r_5 + r_2^2r_4r_8 + 
        r_2r_3^2r_5 + r_2r_3r_5^2 + r_2r_4^2r_8 + r_2r_4r_8^2 + r_3^2r_6r_8 + 
        r_3r_6^2r_8 + r_3r_6r_8^2 + r_4^2r_6r_7 + r_4r_6^2r_7 + r_4r_6r_7^2 + 
        r_5^2r_7r_8 + r_5r_7^2r_8 + r_5r_7r_8^2$,
        
    $\alpha_2=r_1^2r_2r_3 + r_1^2r_4r_6 + r_1^2r_5r_7 + r_1r_2^2r_3 + r_1r_2r_3^2 + r_1r_4^2r_6 
        + r_1r_4r_6^2 + r_1r_5^2r_7 + r_1r_5r_7^2 + r_2^2r_4r_7 + r_2^2r_5r_8 + 
        r_2r_4^2r_7 + r_2r_4r_7^2 + r_2r_5^2r_8 + r_2r_5r_8^2 + r_3^2r_4r_8 + 
        r_3^2r_5r_6 + r_3r_4^2r_8 + r_3r_4r_8^2 + r_3r_5^2r_6 + r_3r_5r_6^2 + 
        r_6^2r_7r_8 + r_6r_7^2r_8 + r_6r_7r_8^2$,
        
    $\alpha_3=r_1^2r_2r_5 + r_1^2r_3r_6 + r_1^2r_4r_7 + r_1r_2^2r_5 + r_1r_2r_5^2 + r_1r_3^2r_6 
        + r_1r_3r_6^2 + r_1r_4^2r_7 + r_1r_4r_7^2 + r_2^2r_3r_4 + r_2^2r_7r_8 + 
        r_2r_3^2r_4 + r_2r_3r_4^2 + r_2r_7^2r_8 + r_2r_7r_8^2 + r_3^2r_5r_8 + 
        r_3r_5^2r_8 + r_3r_5r_8^2 + r_4^2r_6r_8 + r_4r_6^2r_8 + r_4r_6r_8^2 + 
        r_5^2r_6r_7 + r_5r_6^2r_7 + r_5r_6r_7^2$,
        
    $\alpha_4=r_1^2r_2r_4 + r_1^2r_3r_5 + r_1^2r_6r_7 + r_1r_2^2r_4 + r_1r_2r_4^2 + r_1r_3^2r_5 
        + r_1r_3r_5^2 + r_1r_6^2r_7 + r_1r_6r_7^2 + r_2^2r_3r_8 + r_2^2r_5r_7 + 
        r_2r_3^2r_8 + r_2r_3r_8^2 + r_2r_5^2r_7 + r_2r_5r_7^2 + r_3^2r_4r_6 + 
        r_3r_4^2r_6 + r_3r_4r_6^2 + r_4^2r_7r_8 + r_4r_7^2r_8 + r_4r_7r_8^2 + 
        r_5^2r_6r_8 + r_5r_6^2r_8 + r_5r_6r_8^2$
        
        and 
        
 $\beta_1=r_1^2r_2r_3 + r_1^2r_4r_6 + r_1^2r_5r_7 + r_1r_2^2r_7 + r_1r_2r_4^2 + r_1r_3^2r_6 
        + r_1r_3r_5^2 + r_1r_4r_7^2 + r_1r_5r_6^2 + r_2^2r_3r_5 + r_2^2r_4r_8 + 
        r_2r_3^2r_4 + r_2r_5^2r_7 + r_2r_5r_8^2 + r_2r_7^2r_8 + r_3^2r_5r_8 + 
        r_3r_4^2r_6 + r_3r_4r_8^2 + r_3r_6^2r_8 + r_4^2r_7r_8 + r_4r_6^2r_7 + 
        r_5^2r_6r_8 + r_5r_6r_7^2 + r_6r_7r_8^2$,
        
    $\beta_2=r_1^2r_2r_7 + r_1^2r_3r_4 + r_1^2r_5r_6 + r_1r_2^2r_3 + r_1r_2r_5^2 + r_1r_3^2r_5 
        + r_1r_4^2r_7 + r_1r_4r_6^2 + r_1r_6r_7^2 + r_2^2r_4r_7 + r_2^2r_5r_8 + 
        r_2r_3^2r_8 + r_2r_3r_4^2 + r_2r_4r_8^2 + r_2r_5r_7^2 + r_3^2r_4r_6 + 
        r_3r_5^2r_8 + r_3r_5r_6^2 + r_3r_6r_8^2 + r_4^2r_6r_8 + r_4r_7^2r_8 + 
        r_5^2r_6r_7 + r_5r_7r_8^2 + r_6^2r_7r_8$, 
        
    $\beta_3=r_1^2r_2r_4 + r_1^2r_3r_5 + r_1^2r_6r_7 + r_1r_2^2r_5 + r_1r_2r_7^2 + r_1r_3^2r_4 
        + r_1r_3r_6^2 + r_1r_4^2r_6 + r_1r_5^2r_7 + r_2^2r_3r_4 + r_2^2r_7r_8 + 
        r_2r_3^2r_5 + r_2r_3r_8^2 + r_2r_4^2r_7 + r_2r_5^2r_8 + r_3^2r_6r_8 + 
        r_3r_4^2r_8 + r_3r_5^2r_6 + r_4r_6^2r_8 + r_4r_6r_7^2 + r_4r_7r_8^2 + 
        r_5r_6^2r_7 + r_5r_6r_8^2 + r_5r_7^2r_8$, 
        
    $\beta_4=r_1^2r_2r_5 + r_1^2r_3r_6 + r_1^2r_4r_7 + r_1r_2^2r_4 + r_1r_2r_3^2 + r_1r_3r_4^2 
        + r_1r_5^2r_6 + r_1r_5r_7^2 + r_1r_6^2r_7 + r_2^2r_3r_8 + r_2^2r_5r_7 + 
        r_2r_3r_5^2 + r_2r_4^2r_8 + r_2r_4r_7^2 + r_2r_7r_8^2 + r_3^2r_4r_8 + 
        r_3^2r_5r_6 + r_3r_4r_6^2 + r_3r_5r_8^2 + r_4^2r_6r_7 + r_4r_6r_8^2 + 
        r_5^2r_7r_8 + r_5r_6^2r_8 + r_6r_7^2r_8$

\subsection{Vieta's formulas}
Let $d=4a^3+27b^2$.

By using Vieta's formulas relating the roots and the coefficients of a polynomial, we can find the coefficients of the resolvent quartics $h_1(x)=(x-\alpha_1)(x-\alpha_2)(x-\alpha_3)(x-\alpha_4)$, $h_2=(x-\beta_1)(x-\beta_2)(x-\beta_3)(x-\beta_4)$, and $h_3=(x-\gamma_1)(x-\gamma_2)(x-\gamma_3)(x-\gamma_4)$. Again, the formulas involve slightly more than the symmetric polynomials, so there is a bit of additional work to do. Let $\omega_f$ be the evaluation homomorphism from $k[x_1, \ldots, x_8] \rightarrow k_f$ given by $x_i \mapsto r_i$ where the $r_i$s are the roots of $f(x)$.

 Let $u_{35}= \displaystyle \sum_{\sigma \in G} \sigma(x_1^7x_2^5x_3^3x_5)$. Let $v_{35}=-1/8 u_{35} \displaystyle|_{\langle x_1,x_2,x_3,x_4,-x_4,-x_2,-x_3,-x_1 \rangle}=x_1^7x_2^5x_3^3x_4 - x_1^7x_2^5x_3x_4^3 - x_1^7x_2^3x_3^5x_4 + x_1^7x_2^3x_3x_4^5 + x_1^7x_2x_3^5x_4^3 - 
    x_1^7x_2x_3^3x_4^5 - x_1^5x_2^7x_3^3x_4 + x_1^5x_2^7x_3x_4^3 + x_1^5x_2^3x_3^7x_4 - x_1^5x_2^3x_3x_4^7 - 
    x_1^5x_2x_3^7x_4^3 + x_1^5x_2x_3^3x_4^7 + x_1^3x_2^7x_3^5x_4 - x_1^3x_2^7x_3x_4^5 - x_1^3x_2^5x_3^7x_4 + 
    x_1^3x_2^5x_3x_4^7 + x_1^3x_2x_3^7x_4^5 - x_1^3x_2x_3^5x_4^7 - x_1x_2^7x_3^5x_4^3 + x_1x_2^7x_3^3x_4^5 + 
    x_1x_2^5x_3^7x_4^3 - x_1x_2^5x_3^3x_4^7 - x_1x_2^3x_3^7x_4^5 + x_1x_2^3x_3^5x_4^7$.
    
This element is a fundamental invariant and so by Proposition 6.3.3 of Adelmann, \cite{Adelmann}, the coefficients of the resolvent polynomials $\mathcal{RP}_{S_8}(u_{35})(Y)$ are symmetric polynomials in the $x_i$s and can therefore be expressed by elementary symmetric polynomials. By proposition 6.3.8 of Adelmann, the polynomial $\mathcal{R}_{S_8,u_{35}}(Y)$ is a polynomial whose coefficients are polynomials in the coefficients of the defining polynomial and $\omega_f(u_{35})$ appears as a simple linear factor. This polynomial has degree 210 and it too computationally time-consuming. Instead, we can study $\mathcal{RP}_{S_4}(v_{35})(Y)$. This polynomial has degree 2. 

    Then, 
    $$v_{35} = \displaystyle \prod_{i=1}^4 x_i  \displaystyle \prod_{i \ne j} (x_i-x_j)  \displaystyle \prod_{i \ne j} (x_i+x_j)$$
    
$$\omega_f(v_{35})=-2^6(4a^3+27b^2)(27bz^3-9a^2z^2-a^3).$$

Evaluating the $x_i$ at $r_i$, we get the following:\\

$h_1=(x-\alpha_1)(x-\alpha_2)(x-\alpha_3)(x-\alpha_4)=x^4-512dx^2 +2^{15}dw^2x  + 2^{16}d(d+w^2(12az-36b)) +  2^{18}d(27bz^3 - 9a^2z^2 - a^3)$ \\

$h_2=(x-\beta_1)(x-\beta_2)(x-\beta_3)(x-\beta_4)=x^4 -512dx^2 +2^{15}dw^2x  + 2^{16}d(d+w^2(12az-36b)) $ \\

 $h_3=(x-\gamma_1)(x-\gamma_2)(x-\gamma_3)(x-\gamma_4)=x^4-8wx^3+6(2az+3b)x^2-d$.\\

Let $\Delta=-16(4a^3+27b^2)$ and let $g=x^4-4\Delta x-12 a \Delta$, the polynomial for the unique $S_4$-quotient of the $\mathrm{GL}_2(\Z/4\Z)$-extension of $\Q$ given by $\Q(E[4])/\Q$ when $\overline{\rho}_{E,4}$ is surjective. Let $k_g$ be the (generically degree 4) field extension of $\Q$ given by adjoining a root of $g$ to $\Q$. Let $k_1$ be the (generically degree 4) field extension of $\Q$ given by adjoining a root of $h_1$ to $\Q$. We will show in the quaternion origami paper that there is an isomorphism between $k_g$ and $k_1$.

\bibliography{adel}
\bibliographystyle{plain}

\end{document}